\magnification=\magstep1
\vbadness=10000
\hbadness=10000
\tolerance=10000

\def\C{{\bf C}}   
\def\Z{{\bf Z}}   
\def\F{{\bf F}}   
\def\Sp{\hbox{\rm Sp}}
\def\GL{\hbox{\rm GL}}
\def\Aut{\hbox{\rm Aut}}
\def\Ma{\hbox{\rm mass}}
\def\O{\hbox{\rm O}}
\def\pii{\pi\hbox{\rm i}}

\proclaim
A Siegel cusp form of degree 12 and weight 12. 
J. reine angew. Math 494 (1998) 141-153.

Richard E. Borcherds,
\footnote{$^*$}{ Supported by a Royal Society
professorship.}

D.P.M.M.S.,

16 Mill Lane,

Cambridge,
CB2 1SB,
England.

 reb@dpmms.cam.ac.uk

\bigskip

E. Freitag,

Universit\"at Heidelberg,

Im Neuenheimer Feld 288

D-69120 Heidelberg, Germany.

freitag@mathi.uni-heidelberg.de
\bigskip

R.  Weissauer,

Universit\"at Mannheim,

D-68131 Mannheim, Germany.

weissauer@math.uni-mannheim.de

\medskip\noindent
It has been conjectured by Witt [Wi] (1941)
and proved later (1967) independently by Igusa [I]
and M.\ Kneser [K] that
the theta series with respect to the two unimodular even
positive definite lattices of rank 16 are linearly dependent in
degree $\le 3$
and linearly independent in degree 4.
In this paper we consider the next
case of the 24 Niemeier lattices of rank 24. The associated
theta series are linearly dependent in degree $\le 11$ and linearly  
independent
in degree 12. The resulting Siegel cusp form of degree 12 and weight
12 is a Hecke eigenform which seems to have
interesting properties. We would like to thank G.\ H\"ohn for
helpful comments and hints.

\proclaim Construction of Siegel cusp forms by theta series.

\noindent
Let $\Lambda$ be an even unimodular positive definite
lattice, i.e.\ a free abelian group equipped with
a positive definite symmetric bilinear form $(x,y)$,
such that $\Lambda$ coincides with its dual and such that
$$Q(x):={1\over 2}(x,x)$$
is integral. By reduction mod $2$ we obtain
a quadratic form
$$q:E:=\Lambda/2\Lambda\longrightarrow\Z/2\Z,
  \qquad q(a+2\Lambda)=Q(a)\hbox{ \rm mod }2.$$
on the $\Z/2\Z$-vector space $E$.
The standard theta series of degree $n$ with respect to
$\Lambda$ is
$$\vartheta_{\Lambda}(Z)=\sum_{g\in \Lambda^n}\exp\pi\hbox{\rm i}
  \sigma(T(g)Z)\qquad(\sigma=\hbox{trace}),$$
where
$$T(g):=\bigl((g_i,g_j)\bigr)_{1\le i,j\le n}\qquad(g=(g_1,\dots,g_n)).$$
The variable $Z$ varies on the Siegel upper half plane
of degree $n$.
This is a modular form with respect to the full Siegel
modular group $\Sp(2n,\Z)$, but is not a cusp form.
The weight is $m/2$ if $m$ denotes the rank of $\Lambda$,
and $m$ is divisible by $8$.
To obtain  a cusp from we modify this definition.
\smallskip
Assume that a function $\epsilon(F)$ is given which
depends on subspaces $F\subset E$.
For $g\in \Lambda^n$ we denote by $F(g)$ the image of
$\Z g_1+\cdots+\Z g_n$ in $E$.
For an arbitrary degree $n$ we define
$$f^{(n)}(Z):=\sum_{g\in\Lambda^n}\epsilon(F(g))\exp{\pii\over2}
  \sigma(T(g)Z)\qquad(\sigma=\hbox{\rm trace}).$$
In general this will not be a modular form with
respect to the full modular group.
\smallskip
To construct a suitable function  $\epsilon(F)$
we use the orthogonal group $\O(E)$ of
the vector space $E$. It consists of all elements from
the general linear group $\GL(E)$ which preserve the quadratic form $q$.
It is a basic fact for our construction that $\O(E)$ admits
a subgroup of index $2$. It is the kernel of the
so-called {\it Dickson invariant}.
We refer to [B] for some details.
To define the Dickson invariant we chose a basis
$e_1,\dots e_{m}$ of $E$ such that $q$ is of the form
$$q\Bigl(\sum_{i=1}^{m} x_ie_i\Bigr)=\sum_{j=1}^{m/2}x_jx_{m/2+j},$$
which is possible because all even unimodular lattices are equivalent
over $\Z/p\Z$ for any natural number $p$.
The orthogonal group $\O(E)$ now appears as a subgroup of the symplectic
group $\Sp(m,\Z/2\Z)$. It consists of all symplectic
matrices $M=\left(A B\atop C D\right)$
such that the diagonals of $A'C$ and $B'D$ are
zero. This is the image of the so-called theta group.
It is easy to check that
$$D:\O(E)\longrightarrow\Z/2\Z,\quad D(M)=\sigma(C'B),$$
is a homomorphism. This is the Dickson invariant.
It is non-trivial because if $a\in E$ is an element with
$q(a)\ne 0$ then the ``transvection''
$x\mapsto x-(a,x)a$ has non-zero Dickson invariant.
\smallskip
A subspace $F\subset E$ is called {\it isotropic\/}
if the restriction of $q$ to $F$ vanishes.
We now consider maximal isotropic subspaces of $E$.
Their dimension
is $m/2$. The orthogonal group $\O(E)$ acts transitively on the
set of these spaces. But under the kernel of the Dickson invariant $D$
there are two orbits. Two spaces $F_1$ and $F_2$ are in the same
orbit if and only if their intersection has even dimension.
We select one of the two orbits and call it the first orbit
and call the other the second orbit.
\smallskip
We now define a special $\epsilon(F)$ as follows.
It is different from $0$ if and only if
$F$ is maximal isotropic. It is $1$ on the first orbit and $-1$
on the second one.
\smallskip
In the following we consider the system of modular forms $f^{(n)}$
constructed by means of this special $\epsilon(F)$.
\smallskip
Our first observation is that the functions $f^{(n)}(Z)$ have period $1$
in all variables and hence admit a Fourier expansion
$$f^{(n)}(Z)=\sum_T a_n(T)\exp(\pii\sigma(TZ)),$$
where $T$ runs over all integral symmetric matrices with even
diagonal.
Our next observation is that
the coefficients $a_n(T)$ are invariant under unimodular
substitutions $T\mapsto U'TU$, where $U\in\GL(n,\Z)$.
Let $L$ be an arbitrary even lattice of rank $n$. The
Gram matrix $T=\bigl((e_i,e_j)\bigr)$ with respect to a
basis of the lattice is determined up to unimodular
equivalence. We can define
$$a(L):=a_n(T).$$
An easy computation gives
$$a(L)=\#\Aut(L)\sum_{M}\epsilon(M/2M),$$
where the sum is over all $M$ such that
\item{1.} $M$ is a $n$-dimensional sublattice of $\Lambda$.
\item{2.} $M$ is isomorphic to $L(2)$.
($L(2)$ denotes the doubled lattice $L$.
It has the same underlying group as $L$ but the
norms $(x,x)$ are doubled.)
\item{3.} $M/2M$ is maximal isotropic in $\Lambda/2\Lambda$.
\smallskip
The group $\Aut(\Lambda)$ acts on the set of all $M$.
It acts also on the subspaces $F\subset E$. We later need to know
that this group preserves the Dickson invariant.
This is the case if $\Aut(\Lambda)$ is contained
in the special orthogonal group.
For this one has to use that
the composition of the natural homomorphism
$\Aut(\Lambda)\to\O(E)$ with $(-1)^D$ is the determinant [B].
\smallskip
In the following we assume that all automorphisms of $\Lambda$
have determinant $+1$.
Otherwise all $f^{(n)}$ vanish. So we have to exclude
all lattices $\Lambda$ which contain a vector of norm $2$.
We can reformulate the formula
for the Fourier coefficients as

\proclaim Lemma 1. The Fourier coefficients $a(L)$ of
the functions $f^{(n)}$ are given by
$$a(L)=\#\Aut(\Lambda)\#\Aut(L)
   \sum_M{\epsilon(F)\over\#\Aut(\Lambda,M)}\qquad(F=M/2M).$$
Here $M$ runs over a set of representatives of $\Aut(\Lambda)$-orbits
of sublattices  of $\Lambda$ which are isomorphic to $L(2)$.
The group $\Aut(\Lambda,M)$ consists of all elements of
$\Aut(\Lambda)$ which preserve $M$ as a set.

\noindent
We want to prove now that $f^{(n)}$ is a modular form with
respect to the full modular group. More precisely $f:=(f^{(n)})$ is a
stable system of Siegel modular forms, i.e.\ $f^{(n)}$ can be obtained
from $f^{(n+1)}$ by applying the Siegel $\Phi$-operator.
It is known that every stable system can be written
in a canonical way as linear combination of the
standard theta functions $\vartheta_L$. This leads us to
the following construction of a linear combination of
standard theta series.
\smallskip
Let $F\subset E$ be a maximal isotropic space. We consider the inverse
image $\pi^{-1}(F)$ of $F$ under the natural projection
$\pi:\Lambda\to E$. The quadratic form $Q/2$ is even and unimodular
on $\pi^{-1}(F)$. In this way we obtain a new $m$-dimensional
even  unimodular lattice $\Lambda_F$. This
is the so-called perestroika of $\Lambda$ with respect to $F$
in the notation of Koch and Venkov [KV].
\smallskip
We need some more notation. Let $\Lambda'$
be an even unimodular positive definite
lattice of dimension $m$. We introduce the mass and the
modified mass by
$$\eqalign{
\Ma(\Lambda')&=\sum_{\Lambda_F\cong  
\Lambda'}{1\over\#\Aut(\Lambda,F)},\cr
\Ma^\epsilon(\Lambda')&=\sum_{\Lambda_F\cong
\Lambda'}{\epsilon(F)\over\#\Aut(\Lambda,F)},\cr}
$$
where $F$ runs over a system of representatives of $\Aut(\Lambda)$-orbits
of maximal isotropic subspaces of $E$ with
perestroika of type $\Lambda'$.
\smallskip
We fix a system $\Lambda_1,\dots,\Lambda_h$ of representatives
of isomorphism classes of such lattices $\Lambda'$ and write
$$\Ma(i):=\Ma(\Lambda_i),\quad\Ma^\epsilon(i):=\Ma^\epsilon(\Lambda_i).$$

\proclaim Theorem 2. We have
$$f=\#\Aut(\Lambda)\sum_{i=1}^h\Ma^\epsilon(i)\vartheta_{\Lambda_i}.$$
In particular the $f^{(n)}$ are modular forms with respect to the  
full modular
group. The forms $f^{(n)}$ vanish for $n<m/2$, and are cusp
forms for $n=m/2$.

\noindent
{\it Proof}. The right hand side of the  equation in theorem 2
can be written as
$\sum_F\epsilon(F)\vartheta_{\Lambda_F}$. So the difference between both
sides is
$$\sum_F\epsilon(F)\sum_{g:\ F(g)\subset F,\ F(g)\ne F}
  \exp{\pii\over2} \sigma(T(g)Z).$$
We have to show that this series vanishes. We even show
that the partial sum for each fixed $g$ vanishes.
This means:
\smallskip\noindent
{\it Let $F'\subset E$ be an isotropic subspace which is
not maximal. Then
$$\sum_{F'\subset F}\epsilon(F)=0.$$}
This follows from the existence of an element
$g\in\O(E)$
which stabilizes $F'$ and
which has non-trivial Dickson invariant. The  existence
of such a $g$ can be proved easily by using the above
normal form of $E$.
This  proves  Theorem 2.
\smallskip
The main problem is whether the cusp form $f^{(m/2)}$ vanishes
identically or not. This depends on the lattice $\Lambda$. In the next
section we show that it does not vanish if $\Lambda$ is the
Leech lattice.

\proclaim The Siegel cusp form is nonzero in
case of the Leech lattice.

\noindent {}From now on we assume
that $\Lambda$ is the Leech lattice.
In this
section we show that
the Siegel form $f^{(12)}$
from the previous section does not vanish in this case. It is
a Siegel cusp form of degree 12
and weight 12.
Actually we will give two proofs. The first one uses
computer calculations and uses the representation
as linear combination of standard theta series. The second
proof
uses the original definition and
is independent of computer calculations.

\proclaim A first proof for the non vanishing of the cusp form  
$f^{(12)}$.

\noindent
By theorem 2 the Siegel modular form $f^{(12)}$
is a linear combination of
degree 12 theta functions of the 24 Niemeier lattices.
We refer to  [CS] for a detailed description of
the Niemeier lattices. If $L$ is a Niemeier lattice
different from the Leech lattice the vectors of norm $2$
generate a sublattice $L_0$ which determines $L$ up to
isomorphism. We use the notation $L=L_0^+$.
Hence
$D_{24}^+$ is (up to isomorphism) the unique Niemeier lattice
which contains the root lattice $D_{24}$. We use the usual
notations [CS] for the root lattices.
We want to compute the modified
mass of $D_{24}^+$.

\proclaim Lemma 3. The group $\Aut(\Lambda)$ acts transitively
on the set of all sublattices of $\Lambda$ which are isomorphic
to $D_{24}(2)$. The same is true for $D_{24}^+(2)$.

\noindent
{\it Proof}. Recall that a frame of $\Lambda$ is a set of 24  
distinct pairs
$\pm v_i$ of norm 8 vectors of $\Lambda$ all congruent mod $2\Lambda$.
(See lecture 3 of chapter 10 of [CS].)
To every frame we may associate a copy of $D_{24}(2)$ in $\Lambda$.
It is generated by the vectors $(\pm v_i\pm v_j)/2$.
We also get an embedding of $D_{24}^+(2)$ into $\Lambda$ because
the glue vector is contained in $\Lambda$.
It is easy to see that this defines  bijections between
frames and sublattices of type $D_{24}(2)$ and between
frames and sublattices  of type $D_{24}^+(2)$.
The group $\Aut(\Lambda)$
permutes the frames transitively. This proves Lemma 3.
\smallskip
The image $F_0$ of an embedded $D_{24}^+(2)$-lattice
is a maximal isotropic subspace of
$\Lambda/2\Lambda$.
All these $F_0$ have the same Dickson sign because of lemma 3.
We normalize the Dickson sign that it is 1 on these $F_0$.
If $F$ is any other maximal isotropic
subspace then
$$\epsilon(F)=(-1)^{\dim_{\F_2}(F\cap F_0)}.$$
{}From lemma 3 now follows that the mass and the modified mass
with respect to $D_{24}^+$ agree.
The 24 masses $\Ma(i)$
have been computed in the paper [DLMN].
The Niemeier lattice $D_{24}^+$ has index $i=24$.
We will denote the number of frames by $n_F=8292375=3^6.5^3.7.13$.
{}From
Lemma 3 we obtain
$$\Ma(24)=\Ma^\epsilon(24)={n_F\over\#\Aut(\Lambda)}=
  {1\over 2}\cdot{1\over 501397585920}$$
in accordance with [DLMN].
\smallskip
The following table shows a $24\times 24$-matrix.
The columns correspond to the 24
Niemeier lattices
in order of their Coxeter numbers.
This is the same order as used in [DLMN] and which we
will use in Theorem 4 below.
The rows correspond to
the following lattices of degree $\le12$:
the 0 dimensional lattice, $A_j$ for $1\le j\le11$,
$D_j$ for $4\le j\le12$, and $E_6,E_7,E_8$.
The matrix entry is the number of sublattices of
the Niemeier lattice
isomorphic to the lattice of each row.
This means that each column contains
suitable normalized
Fourier coefficients
of the theta function of the corresponding Niemeier lattice.
\def\s{\scriptscriptstyle}\bigskip\noindent%
\medskip\noindent
\vbox{\baselineskip=6pt
\halign{$\s#\ $\hfill&
   $\s#\,$\hfill&$\s#\,$\hfill&$\s#\,$\hfill&$\s#$\
 \hfill&$\s#\,$\hfill&$\s#\,$
\hfill&$\s#\,$\hfill&$\s#\,$\hfill&$\s#\,$\hfill&$\s#$\
 \hfill&$\s#\,$\hfill&$\s#\,$
\hfill&$\s#\,$\hfill&$\s#\,$\hfill&$\s#\,$\hfill&$\s#$\
 \hfill&$\s#\,$\hfill&$\s#\,$
\hfill&$\s#\,$\hfill&$\s#\,$\hfill&$\s#\,$\hfill&$\s#$\
 \hfill&$\s#\,$\hfill&$\s#$\cr
0&1&1&1&1&1&1&1&1&1&1&1&1&1&1&1&1&1&1&1&1&1&1&1&1\cr
A_1&0&24&36&48&60&72&72&84&96&108&120&120&144&144&156&168&192&216&216&264&300&360&360&552\cr  
A_2&
0&0&12&32&60&96&96&140&192&252&320&320&480&480&572&672&896&1152&1152&1760&2300&3360&3360&8096\cr  
A_3&
0&0&0&8&30&72&72&140&240&378&560&560&1080&1080&1430&1848&2912&4320&4320&8360&12650&22680&22680&87032\cr  
A_4&
0&0&0&0&6&24&0&84&144&378&600&384&1344&864&2574&2688&6384&8064&10584&25344&53130&94080&72576&680064\cr  
A_5&
0&0&0&0&0&4&0&28&56&252&452&128&1184&144&3432&2688&10696&9408&19908&59136&177100&296576&120960&4307072\cr  
A_6&
0&0&0&0&0&0&0&4&16&108&240&0&856&0&3432&1536&13744&8256&32112&101376&480700&766720&103680&22150656\cr  
A_7&
0&0&0&0&0&0&0&0&2&27&90&0&495&0&2574&384&14022&5832&43794&126720&1081575&1660320&38880&94140288\cr  
A_8&
0&0&0&0&0&0&0&0&0&3&20&0&220&0&1430&0&11696&2560&48620&112640&2042975&2929600&2880&334721024\cr  
A_9&
0&0&0&0&0&0&0&0&0&0&2&0&66&0&572&0&8008&512&43758&67584&3268760&4100096&0&1004163072\cr  
A_{10}&
0&0&0&0&0&0&0&0&0&0&0&0&12&0&156&0&4368&0&31824&24576&4457400&4472832&0&2556051456\cr  
A_{11}&
0&0&0&0&0&0&0&0&0&0&0&0&1&0&26&0&1820&0&18564&4096&5200300&3727360&0&5538111488\cr  
D_4&
0&0&0&0&0&1&6&0&10&0&15&60&80&180&0&210&126&840&315&990&0&4970&9450&10626\cr  
D_5&
0&0&0&0&0&0&0&0&2&0&6&24&48&108&0&168&126&1008&378&1584&0&11928&22680&42504\cr  
D_6&
0&0&0&0&0&0&0&0&0&0&1&4&7&0&0&84&84&336&63&1848&0&11788&11340&134596\cr  
D_7&
0&0&0&0&0&0&0&0&0&0&0&0&1&0&0&24&36
 &120&0&1584&0&12520&3240&346104\cr D_8&
0&0&0&0&0&0&0&0&0&0&0&0&0&0&0&3&9&45&0&990&0&13005&405&735471\cr
D_9&
0&0&0&0&0&0&0&0&0&0&0&0&0&0&0&0&1&10&0&440&0&11440&0&1307504\cr
D_{10}&
0&0&0&0&0&0&0&0&0&0&0&0&0&0&0&0&0&1&0&132&0&8008&0&1961256\cr
D_{11}&
0&0&0&0&0&0&0&0&0&0&0&0&0&0&0&0&0&0&0&24&0&4368&0&2496144\cr
D_{12}&
0&0&0&0&0&0&0&0&0&0&0&0&0&0&0&0&0&0&0&2&0&1820&0&2704156\cr E_6&
0&0&0&0&0&0&0&0&0&0&0&0&1&4&0&0&0&56&28&0&0&1120&3360&0\cr E_7&
0&0&0&0&0&0&0&0&0&0&0&0&0&0&0&0&0&2&1&0&0&120&360&0\cr E_8&
0&0&0&0&0&0&0&0&0&0&0&0&0&0&0&0&0&0&0&0&0&1&3&0\cr}}
\bigskip\noindent
The rank of this matrix is $24$.
Therefore the $24$ theta functions are linearly independent in
degree $12$.
All lattices of rows in the above table other than $D_{12}$
are of rank $<12$. {}From this we
see  that
the space of degree 12 cusp forms
spanned by the 24 theta functions is at most one dimensional.
\smallskip
In the following we use the notation
$\vartheta(L_0)=\vartheta_L$
for a root lattice $L_0$ contained in a Niemeier lattice $L=L_0^+$.
We deduce from the above matrix:
\proclaim Theorem 4. The $24$ theta functions are
linearly independent in degree $12$.
Every degree $12$ cusp form
spanned by the $24$ theta functions is a constant multiple of
\medskip\noindent
\halign{\vphantom{\Big(}\  
\hfil$#$&$#$\hfil&\hfil$#$&$#$\hfil&\hfil$#$&$#$\hfil\cr
 {1\over 152769576960} &\vartheta(\hbox{\rm Leech})
&-{1\over 3183476736} &\vartheta(A_1^{24})
&+{1\over 591224832} &\vartheta(A_2^{12})\cr
-{5\over 1146617856} &\vartheta(A_3^8)
&+{13\over 1990656000} &\vartheta(A_4^6)
&-{83\over 11943936000} &\vartheta(A_5^4 D_4)\cr
+{19\over 16307453952} &\vartheta(D_4^6)
&+{41\over 15676416000} &\vartheta(A_6^4)
&-{1\over 37158912000} &\vartheta(A_7^2 D_5^2)\cr
-{197\over 351151718400} &\vartheta(A_8^3)
&+{59\over 214990848000} &\vartheta(A_9^2 D_6)
&-{13\over 229323571200} &\vartheta(D_6^4)\cr
-{1\over 35831808000} &\vartheta(A_{11} D_7 E_6)
&+{1\over 143327232000} &\vartheta(E_6^4)
&+{31\over 7685922816000} &\vartheta(A_{12}^2)\cr
+{37\over 11415217766400} &\vartheta(D_8^3)
&-{29\over 21069103104000} &\vartheta(A_{15} D_9)
&-{1\over 7023034368000} &\vartheta(D_{10} E_7^2)\cr
+{1\over 3511517184000} &\vartheta(A_{17} E_7)
&+{53\over 4237899595776000} &\vartheta(D_{12}^2)
&-{1\over 1332620771328000} &\vartheta(A_{24})\cr
-{1\over 3595793596416000} &\vartheta(D_{16} E_8)
&+{1\over 10787380789248000} &\vartheta(E_8^3)
&+{1\over 2729207339679744000} &\vartheta(D_{24})\cr}
\bigskip\noindent
{\bf Remark.} This
linear combination has been normalized
so that the Fourier coefficient
of the $D_{12}$ lattice is 1.

\noindent
We know that at least one coefficient in theorem 2 is different
from $0$.
{}From the linear independence of theorem 4 we obtain:

\proclaim Theorem 5.
The form $f^{(12)}$ of theorem 2 does not vanish when $\Lambda$ is the
Leech lattice.

\noindent
We know already that $f^{(12)}$ is a cusp form. This cusp form
must be a constant multiple of the linear combination of theta  
functions in
theorem 4. The constant factor can be determined by looking at
the coefficient of $D^+_{24}$:

\proclaim Theorem 6. Let $g$ be the linear combination
from theorem 4. This linear combination is
a cusp form. We have
$$f=\#\Aut(\Lambda)\cdot2^6\cdot3^5\cdot5^2\cdot7\cdot g.$$

\noindent
The following corollary is due to Igusa, who gave in [I] an affirmative
answer to a question asked by Witt. Igusa used deep results about
modular forms. An elementary proof has been found a litle later
by M.\ Kneser [K].

\proclaim Corollary. If $E_8^2$ and $D_{16}^+$ are the two
16 dimensional even unimodular lattices and $L$ is a lattice
of dimension at most 3 then the numbers of embeddings of
$L$ into $E_8^2$ and $D_{16}^+$ are the same.

\noindent
{\it Proof.} Look at the Fourier coefficient in $f$ of the lattice
$M=L\oplus E_8$. This coefficient is given in terms of
the numbers of embeddings of $M$ into various Niemeier lattices,
and must vanish as $f$ is a cusp form of
degree 12 and $M$ has dimension at most 11. On the other hand any
Niemeier lattice containing $M$  must contain an $E_8$ sublattice so
must be $E_8^3$ or $D_{16}^+E_8$, and the number of embeddings of $M$
into $E_8^3$ or $D_{16}E_8$ is given up to some fixed factors by the
number of embeddings of $L$ into $E_8^2$ or $D_{16}^+$
respectively. This easily implies that the numbers of embeddings of
$L$ into $E_8^2$ or $D_{16}^+$ are equal.
This proves the  corollary.

\proclaim A second proof for the non vanishing of the cusp form.

\noindent
We will give another proof
that $f^{(12)}$ is not identically 0, without using
computer calculations,
by showing
that $a(M)$ is nonzero if $M$ is the $D_{12}$ lattice.
We use the formula for the Fourier coefficients from lemma 1
coming from the first representation of $f$.

\proclaim Lemma 7. There is a maximal isotropic subspace $F_0$
of $\Lambda/2\Lambda$ which has an odd number of conjugates under
$Aut(\Lambda)$.

\noindent
{\it Proof.} We take $F_0 $ to be the subspace of $\Lambda/2\Lambda$
spanned by the vectors of the form $\sum_i4n_ie_i$
where the sum is over the usual orthogonal basis $e_i$ of norm
$1/8$ vectors of $\Lambda\otimes R$
[CS], p.287
and the $n_i$'s are either
all in $\Z$ or all in $\Z+1/2$ and have even sum.
Then $F_0$ is fixed by the standard subgroup $2^{12}.M_{24}$
of $Aut(\Lambda)$ which has odd index, so
that $F_0$ has an odd number
of
conjugates under $Aut(\Lambda)$. This proves lemma 7.
\smallskip
The lattice $D_{12}$ is the 12 dimensional  lattice
of determinant 4 generated by the roots of as $D_{12}$ root system,
and  $D_{12}(2)$ is this lattice with all
inner products
multiplied by 2.
\smallskip
As we have mentioned,
a frame of $\Lambda$ is a set of 24 distinct pairs
$\pm v_i$ of norm 8 vectors of $\Lambda$ all congruent mod $2\Lambda$.
For each frame there is an  action of the Mathieu group $M_{24}$
on this 24 element set.

\proclaim Lemma 8. Sublattices $L$ of $\Lambda$ isomorphic to
$D_{12}(2)$
correspond to frames in $\Lambda$ together with
a choice of 12 element subset of the 24 element subset of the frame. 

\noindent
{\it Proof.} Twice the images of the norm 1
vectors of the dual of the $D_{12}$
in $\Lambda$ give a set of 12 pairs $\pm v_i$ of norm 8 vectors of  
$\Lambda$,
all congruent mod 2, and these determine a unique frame
and a 12 element subset of the 24 element set of the frame.
The image $L$ of the $D_{12}(2)$ is then spanned by the set of vectors of
the form $(\pm v_i\pm v_j)/2$.

Conversely any choice of frame and 12 element subset gives
a sublattice isomorphic to $D_{12}(2)$ by the construction above.
This gives a one to one correspondence between
such sublattices and pairs consisting of a frame and a 12 element
subset, which proves lemma 8.

\proclaim Lemma 9. There are exactly 5 orbits of sublattices
of $\Lambda$ isomorphic to $D_{12}(2)$, of sizes
$2^2.5.7.11.23n_F$, $2^4.3^2.5.7.11.23n_F$, $2^6.3^2.7.11.23n_F$,
$2^8.3^2.7.23n_F$, and $2^4.7.23n_F$.

\noindent
{\it Proof.} This follows immediately from lemma 3 and the  
classification of
the five $M_{24}$-orbits of 12 element subsets in [CS, chapter 10,
theorem 22] and the fact that $Aut(\Lambda)$ acts transitively on the
$n_F$ frames.

\proclaim Lemma 10. A sublattice $L$ of $\Lambda$ isomorphic to
$D_{12}(2)$ represents a maximal isotropic subspace of
$\Lambda/2\Lambda$ if and only if there are no vectors $v\in \Lambda$
with $2v\in L$, $v\notin L$.

\noindent
{\it Proof.} The image of $L$ in $\Lambda/2\Lambda$ is isotropic
as all vectors of $L$ have norm divisible by 4,
so we have to check whether this image has dimension 12, in other words
we have to check whether the map from $L/2L$ to $\Lambda/2\Lambda$
is injective. But this is the same as asking whether there
exist no vectors $v$ as above, and this proves lemma 10.

\proclaim Lemma 11. There is a vector $v\in \Lambda$
with $2v\in L$, $v\notin L$, if and only if
the 12 element set $S$ corresponding to $L$ contains no
nonzero elements of the Golay code.

\noindent
{\it Proof.} If there is such a set $S$, then we can construct $v$ as
a sum $\sum_{i\in S}\pm v_i/4$. Conversely given $v$, we can
construct
$S$ as the set of coordinates where the coefficient of $v$
is $2\bmod 4$ (which is nonempty as $v\notin L$). This proves
lemma 11.

\proclaim Lemma 12. There are exactly two $Aut(\Lambda)$-orbits of  
sublattices
$L$ of $\Lambda$ isomorphic to $D_{12}(2)$ such that
$L/2L$ is a maximal isotropic subspace of $\Lambda/2\Lambda$,
and they have  sizes
$2^6.3^2.7.11.23n_F$ and $2^8.3^2.7.23n_F$.

\noindent
{\it Proof.} By lemmas 4, 5 and 6 we have to find the orbits of 12  
element
subsets of a 24 element set acted on by $M_{24}$ which contain no
nonzero elements of the Golay code. These can be read off from the
discussion in section 2.6 of chapter 10 of [CS]. In Conway's
terminology, the ``special'' and ``extraspecial'' 12 element sets (and
no others) contain octads of the Golay code, and the ``umbral'' 12
element sets are already in the Golay code. This leaves the
``transverse'' and ``penumbral'' as the 12 element sets
containing no element of the Golay code, and they have the orbit
sizes as stated in the lemma.  This proves lemma 12.
We now obtain a new proof for

\proclaim Theorem.
In case of the Leech lattice
the Siegel cusp form $f^{(12)}$ of theorem 2 is
nonzero.

\noindent
{\it Proof.}
We use the formula of lemma 1 to compute the Fourier coefficient
$a(D_{12})$. The Dickson signs
are constant on each of the two orbits.  Using the formulas
$$\#\Aut(D_{12})=2^{12}\cdot 12!,\qquad
\#\Aut(\Lambda)=
  2^{22}\cdot3^9\cdot5^4\cdot7^2\cdot11\cdot13\cdot23$$
we obtain
$$a(D_{12})=2^{28}\cdot 3^{13}\cdot 5^5\cdot 7^3\cdot11\cdot13
\cdot23\cdot(\pm 11\pm 4)$$
Whatever the signs of the Dickson invariants might be, this
number is different from $0$. This proves the theorem.
({}From Lemma 3 and 8 it follows  that in fact both signs
are $+1$. This is in accordance with the formula
$a(D_{12})=2^{28}\cdot 3^{14}\cdot 5^6\cdot 7^3\cdot11\cdot13
\cdot23$, which comes from lemma 1.)
\smallskip
The coefficient we have shown is nonzero is in some sense
the first
nonzero coefficient, or more precisely a nonzero coefficient
corresponding to a 12 dimensional lattice of smallest possible
determinant.
This follows easily from the fact that there are no
even 12 dimensional lattices with determinant less than 4.
\smallskip
We observed that the Fourier coefficient of a lattice $M$ of determinant
$n$ often seems to be closely related to the coefficient of $q^n$ of the
weight $13/2$ ordinary modular form
$$\eqalign{
\eta(8 \tau)^{12}  \theta(\tau) =& q^4 + 2q^5 +
2q^{8} - 12q^{12} - 22q^{13} - 24q^{16} + 56q^{20} + 84q^{21} +  
108q^{24}
\cr&-
112q^{28} - 66q^{29} - 176q^{32} + 9q^{36} - 398q^{37} - 196q^{40}
\cr&+
364q^{44} + 990q^{45} + 1056q^{48} - 616q^{52} + 70q^{53} - 728q^{56} 
\cr&+
432q^{60} - 2354q^{61} - 1472q^{64} - 240q^{68} + 1080q^{69} + 990q^{72} 
\cr&-
484q^{76} + 1848q^{77} + 2752q^{80} + 2352q^{84} + 2292q^{85} +  
1276q^{88}
\cr&-
2608q^{92} - 3852q^{93} - 9504q^{96} + O(q^{100})\cr
}$$
at least when $n$ is 0, 4, or 5 mod 8.
It is often the same when $n$ is divisible by 4 and often differs  
by a factor
of $-2$ when $n\equiv 5\bmod 8$. We have not been able
to find such similar properties when $n\equiv 1\bmod 8$.
\smallskip
Here is a table of the Fourier coefficients corresponding
to lattices $M$ that have determinant at most 96 and are generated
by their norm 2 vectors. This table was calculated using the expression
of $f$ as a linear combination of theta functions, and using the
fact that if $M$ is generated by norm 2 vectors then the number
of embeddings of $M$ into any other lattice $L$ can easily be
worked out knowing the root system of $L$. A larger version of this
table and the programs used to calculate it can be found on
R. E. Borcherds' home page http://www.dpmms.cam.ac.uk/\~{}reb.
\medskip\noindent
    \halign{\hfil#&\ \hfil$#$&\ \hfil$#$\ &
\qquad\qquad\hfil#&\ \hfil$#$&\ \hfil$#$\ &
\qquad\qquad\hfil#&\ \hfil$#$&\ \hfil$#$\cr
  det&        \hbox{ coef} &    \hbox{ Lattice} &
  det&        \hbox{ coef} &    \hbox{ Lattice} &
  det&        \hbox{ coef} &    \hbox{ Lattice} \cr
    4&                   1&      D_{12} &
    4&                   1&      D_4 E_8 &
    5&                  -1&      A_4 E_8 \cr
    8&                   2&      A_1 A_3 E_8 &
    8&                   2&      A_1 D_{11} &
    8&                   2&      D_5 E_7 \cr
    9&                   6&      A_2^2 E_8 &
    9&                   6&      E_6^2 &
   12&                 -12&      A_1^2 A_2 E_8 \cr
   12&                 -12&      A_2 D_{10} &
   12&                 -12&      A_5 E_7 &
   12&                 -12&      D_6 E_6 \cr
   13&                  11&      A_{12} &
   16&                  40&      A_1 D_4 E_7 &
   16&                  40&      A_1^2 D_{10} \cr
   16&                  40&      A_1^4 E_8 &
   16&                  40&      D_6^2 &
   16&                 -24&      A_3 D_9 \cr
   16&                 -24&      D_5 D_7 &
   16&                 -88&      D_4 D_8 &
   20&                  -8&      A_1 A_4 E_7 \cr
   20&                  56&      A_4 D_8 &
   21&                 -42&      A_6 E_6 &
   24&                 108&      A_1 A_{11} \cr
   24&                 108&      A_1 A_2 D_9 &
   24&                 108&      A_1 D_5 E_6 &
   24&                 108&      A_2 A_3 E_7 \cr
   24&                 108&      A_5 D_7 &
   28&                -112&      A_6 D_6 &
   32&                 -48&      A_1 A_3 D_8 \cr
   32&                 -48&      A_1 D_4 D_7 &
   32&                 -48&      A_1 D_5 D_6 &
   32&                 -48&      A_1^2 A_3 E_7 \cr
   32&                 -48&      A_1^3 D_9 &
   32&                -176&      A_7 D_5 &
   33&                  54&      A_2 A_{10} \cr
   36&                   9&      A_8 D_4 &
   36&                  48&      A_2^2 D_8 &
   36&                 738&      A_2 D_4 E_6 \cr
   36&                -336&      A_1 A_2^2 E_7 &
   36&                -336&      A_1 A_5 E_6 &
   40&                -196&      A_1 A_4 D_7 \cr
   40&                -196&      A_3 A_9 &
   44&                 364&      A_1^2 A_{10} &
   45&                 234&      A_2 A_4 E_6 \cr
   45&                -495&      A_4 A_8 &
   48&                 288&      A_1 A_5 D_6 &
   48&                 288&      A_1^2 A_2 D_8\cr
   48&                 288&      A_1^2 D_4 E_6 &
   48&                 288&      A_1^3 A_2 E_7 &
   48&                 288&      A_2 D_4 D_6 \cr
   48&                 288&      A_5 A_7 &
   48&                -480&      A_2 A_3 D_7 &
   48&                1056&      A_2 D_5^2 \cr
   48&                1056&      A_3^2 E_6 &
   49&                1260&      A_6^2 &
   56&                -728&      A_1 A_6 D_5 \cr
   60&                 432&      A_1 A_2 A_9&
   60&                 432&      A_1^2 A_4 E_6 &
   60&                 432&      A_2 A_4 D_6 \cr
   64&                1088&      A_1^2 A_3 D_7 &
   64&                1088&      A_1^2 D_4 D_6 &
   64&                1088&      A_1^2 D_5^2 \cr
   64&                1088&      A_1^4 D_8 &
   64&                1088&      A_1^5 E_7 &
   64&                1088&      A_3^2 D_6 \cr
   64&                4160&      A_3 D_4 D_5 &
   64&               -3008&      A_1 A_7 D_4 &
   64&               23616&      D_4^3 \cr
   72&                 990&      A_1 A_3 A_8 &
   72&                -468&      A_1 A_2 A_3 E_6 &
   72&                -468&      A_2 A_5 D_5 \cr
   72&                2448&      A_1 A_2^2 D_7 &
   80&                2240&      A_1 A_4 A_7 &
   80&                2752&      A_3 A_4 D_5 \cr
   80&                6336&      A_1^3 A_9 &
   80&               -1856&      A_1^2 A_4 D_6 &
   80&               -9024&      A_4 D_4^2 \cr
   81&                5886&      A_2^2 A_8 &
   81&               -7236&      A_2^3 E_6 &
   84&                -336&      A_1 A_5 A_6 \cr
   84&               -3024&      A_2 A_6 D_4 &
   96&                4320&      A_2 A_3 A_7 &
   96&               -2592&      A_1 A_2 A_3 D_6 \cr
   96&               -2592&      A_1 A_2 D_4 D_5 &
   96&               -2592&      A_1^2 A_5 D_5 &
   96&               -2592&      A_1^3 A_2 D_7 \cr
   96&               -2592&      A_1^3 A_3 E_6 &
   96&               -2592&      A_3 A_5 D_4 &
&&\cr
}
\bigskip

\proclaim A Hecke eigenvalue.

\noindent
Our cusp form is an eigenform of all Hecke operators because
Hecke operators map cusp forms to cusp forms and preserve
the space generated by standard theta series.
\smallskip
We use the definition of the Hecke operator $T(p)$
as given in [F]. The explicit formula of the action
of $T(p)$ for a prime $p$  on theta series ([F, IV.5.10])  states that
$$\vartheta_{\Lambda_\nu}\vert T(p)=
  \beta(p,m,n)\cdot\sum_{\mu=1}^h n(\Lambda_\mu(p),\Lambda_\nu)
 \vartheta_{\Lambda_\mu},$$
where $n(L(p),M)$ denotes the number of sublattices
of $M$ of type $L(p)$.
The constants $\beta(p,m,n)$ depend on the normalization
of $T(p)$. We refer to [F] for explicit expressions.
For example
$$\beta(p,24,12)=p^{{n(n+1)\over2}-12n}.$$
It is well known and easy to
prove that the matrix with entries
$\#\Aut(\Lambda_\mu)\cdot n(\Lambda_\nu,\Lambda_\mu(p))$
is symmetric.
\smallskip
Because of the linear independence of the 24 theta series
it is sufficient to know one row (or column)
of this matrix if one wants
to compute the eigenvalue $\lambda(p)$ of $T(p)$.
We can compute
one row of this matrix
in the case $p=2$. One easily
derives
$$\Ma(L)={n(L(2),\Lambda)\over\#\Aut(\Lambda)},$$
where $\Ma(L)$ is the mass introduced above. As we already mentioned
these masses have been computed in [DLMN] for the Leech lattice
$\Lambda$.
Using table I of this paper one obtains
$$\lambda(2)=2^7\cdot 3^{11}\cdot5\cdot 17\cdot 901141\cdot
  \beta(2,24,12).$$
We obtain now some information about the Satake parameters
$x_1,\dots x_{12}$ of our cusp form at the place $2$.
We recall briefly their definition [A]. The local Hecke algebra
at a prime is isomorphic to the  ring of invariants
$\C[X_0^{\pm1},\dots,X_n^{\pm1}]^{W_n}$, where $W_n$ is the
symplectic Weyl group ([F], IV.3.19). Every homomorphism
of this ring into the field of complex numbers is the restriction
of  a homomorphism of the whole $\C[X_0^{\pm1},\dots,X_n^{\pm1}]$.
The images $x_j$ of the variables $X_j$ are the Satake parameters.
They are determined up to the action of $W_n$. Every Siegel
eigenform of the local Hecke algebra at a prime $p$
defines such a homomorphism. The Ramanujan conjecture
says that the Satake parameters $x_1,\dots,x_n$ of an eigen cusp  
form have
absolute value $1$. It is known that in degree $n>1$
this is not always true.
\smallskip
To obtain information about the Satake parameters
(at the prime $p=2$) we need the image of the operator $T(p)$
in the local Hecke algebra. This formula can be found in
[F]. We
choose a root $y_j=\sqrt{ x_j}$ for
each Satake parameter. A direct consequence of
formula [F], IV.3.14, a) and b) is
$${\lambda(p)^2\over x_0^{-2}x_1\cdots x_n}=\prod_{j=1}^{12}
\bigl(y_j+y_j^{-1}\bigr)^2\quad\hbox{\rm and}\quad
p^{{n(n+1)\over2}-12n}=x_0^{-2}x_1\cdots x_n.$$
The computed value $\lambda(2)$ now gives:

\proclaim Theorem 13. The Satake parameters $x_i=y_i^2$ of our
cusp form of degree $12$ and weight $12$
at the place $p=2$ satisfy
$$\left\vert\prod_{i=1}^{12}(y_i+y_i^{-1})\right\vert=
{3^{11}\cdot5\cdot17\cdot901141\over 2^{26}}.$$
\smallskip\noindent
{\bf Corollary.} The Ramanujan conjecture $\vert x_i\vert=1$
is violated for $p=2$.

\proclaim Open problems.

\noindent
We list a few questions about the Siegel cusp form $f$.
\smallskip
\item{1.} Are the coefficients of
the cusp form of weight 12 and degree 12 all integers when
normalized so that the coefficient of $D_{12}$ is 1?
\smallskip\item{}
One can prove that the coefficients of $f^{(m/2)}/\#\Aut(\Lambda)$
are contained in $\Z[1/2]$ for arbitrary $\Lambda$ and
in $(1/2)\Z$ in case of the Leech lattice. This means
that the denominators of the normalized coefficients divide
$2^7\cdot3^5\cdot5^2\cdot7$.
\smallskip
\item{2.} Why are the coefficients of $f$ similar to those
of the modular form above? Is there a similar relation for
the coefficients of lattices of determinant $1\bmod 8$?
Is it possible to write down some simple explicit
formula for the coefficients of $f$?
\smallskip\item{}
{}From  [We] it follows that the standard $L$-function $L(f,s)$
of $f$
has a pole at $s=1$. This suggests that $L(f,s)=\zeta(s)L(s)$,
where $L$ belongs to a $24$-dimensional $l$-adic
Galois representation. This Galois
representation cannot be pure (theorem 13) and therefore
one might expect that its weight filtration
sheds light on the relationship with
$\eta(8 \tau)^{12}  \theta(\tau)$.
\bigskip\noindent
\proclaim References.

\item{[A]} A.\ Andrianov.
On zeta-functions of Rankin type associated with Siegel
modular forms,
Modular Functions of One Variable VI, 325-338,
Lecture Notes in Mathematics {\bf 627}, Springer-Verlag
Berlin Heidelberg New York, 1977
\smallskip
\item{[B]} N.\ Bourbaki. \'El\'ements
de Math\'ematique,
Livre VI, Deuxi\'eme Partie, Groupes et Alg\'ebre de Lie,
1.\  Hermann, Actualit\'es Scientifique
et Industrielles, {\bf 1293}, 1962.
\smallskip
\item{[CS]} J.\ H.\ Conway, N.\ J.\ A.\ Sloane.
Sphere packings, lattices
and groups. Springer-Verlag New York,
Grundlehren der mathematischen
Wissenschaften {\bf 290}, 1988
\smallskip
\item{[DLMN]}
C.\ Dong, H.\ Li, G.\ Mason, S.\ P.\ Norton.
Associative subalgebras of the Griess algebra and related topics,
preprint, http://xxx.lanl.gov/abs/q-alg/9607013, 1996
\smallskip
\item{[F]} E.\ Freitag. Siegelsche Modulformen,
Springer-Verlag Berlin New York,
Grundlehren der
mathematischen Wissenschaften {\bf 254}, 1983
\smallskip
\item{[I]} J-I.\ Igusa.\ Modular forms and projective invariants,
Am.\ J.\ of Math. {\bf 89}, 817--855, 1967.
\smallskip
\item{[K]} M.\ Kneser. Lineare Relationen zwischen Darstellungsanzahlen
quadratischer Formen, Math.\ Annalen {\bf 168}, 31--36, 1967
\smallskip
\item{[KV]} H.\ Koch, B.B\ Venkov.\ \"Uber ganzzahlige unimodulare
euklidische Gitter,
J.\ Reine Angew.\ Math. {\bf 398}, 144--168, 1989
\smallskip
\item{[We]} R.\ Weissauer.
Stabile Modulformen und Eisensteinreihen,
Lecture notes  in Mathematics {\bf 1219},
Springer-Verlag Berlin Heidelberg New York,
1986
\smallskip
\item{[Wi]} E.\ Witt. Eine Identit\"at zwischen
Modulformen zweiten Grades, Math.\ Sem.\ Hamburg
{\bf 14}, 323--337, 1941
\bye